\documentclass[12pt]{amsart}
\usepackage{epsf,righttag,latexsym,amssymb,amscd,xspace}
\newcommand{\titlestring}{Partial Resolutions of Orbifold
  Singularities via Moduli Spaces of HYM-type Bundles}
\newcommand{\runningheadstring}{Orbifold Singularities and Moduli of
  Bundles} \newcommand{\msc}{(1991): 32S45 (Primary) 32L07, 53C07,
  32M05 (Secondary)}
\newcommand{\name}{Alexander V.\ Sardo Infirri}
\newtheorem{thm}{Theorem}[section]
\newtheorem{nonumberthm}{Theorem}  
  
\newtheorem{prop}[thm]{Proposition}
\newtheorem{cor}[thm]{Corollary}
\newtheorem{lemma}[thm]{Lemma}

 \theoremstyle{definition}

\newtheorem{conj}{Conjecture}
\newtheorem{nonumberconj}{Conjecture}  
\newtheorem{example}[thm]{Example}

\theoremstyle{remark}

\newtheorem{rmk}[thm]{Remark}      %   \renewcommand{\thermk}{}

\newtheorem{note}{Note}

\newcommand{\C}{{\mathbb C}}
\newcommand{\CX}{{\C^*}}
\newcommand{\R}{{\mathbb R}}
\newcommand{\Z}{{\mathbb Z}} 
\newcommand{\Q}{{\mathbb Q}} 
\newcommand{\N}{{\mathbb N}}

\newcommand{\PP}{{\mathbb P}}

\newcommand{\cA}{{\mathcal A}}

\newcommand{\cG}{{\mathcal G}}
\newcommand{\cH}{{\mathcal H}}

\newcommand{\cN}{{\mathcal N}}
\newcommand{\cO}{{\mathcal O}}

\newcommand{\bgamma}{\boldsymbol{\gamma}}
\newcommand{\bomega}{\boldsymbol{\omega}}

\newcommand{\bdg}{{\mathbf g}}
\newcommand{\bdh}{{\mathbf h}} % obsolete

\newcommand{\bdC}{{\mathbf C}}

\newcommand{\bdK}{{\mathbf K}}
\newcommand{\bdQ}{{\mathbf Q}}

\newcommand{\lieu}{{\mathfrak u}}

\newcommand{\Mga}{M^{\Gamma }}
\newcommand{\glg}{\GL^\ga\!}
\newcommand{\ug}{\U^\ga\!}

\newcommand{\cnga}{\cN^{\Gamma }}
\newcommand{\cngamuo}{\cN^{\Gamma }\cap\mu^{-1}(0)}

\newcommand{\clos}{\overline}

\newcommand{\ga}{\Gamma}
\newcommand{\gat}{\mbox{$\Gamma$}}

\newcommand{\dbar}{\bar{\partial}}

\newcommand{\qbar}{{\bar{q}}}
\newcommand{\zbar}{{\bar{z}}}
\newcommand{\ibar}{{\bar{\imath}}}
\newcommand{\jbar}{{\bar{\jmath}}}
\newcommand{\kbar}{{\bar{k}}}

\newcommand{\gitquot}[1]{/\!\!/_{\!#1}\,}

\newcommand{\ie}{i.e.\xspace}
\newcommand{\cf}{c.f.\xspace}
\newcommand{\resp}{resp.\xspace}

\newcommand{\nbd}{neighbourhood}

\newcommand{\git}{geometric invariant theory}

\newcommand{\hip}{positive definite Hermitian inner product}

\newcommand{\kah}{K\"ahler}

\newcommand{\hym}{Hermitian-Yang-Mills}

\DeclareMathOperator{\Ind}{Ind}

\DeclareMathOperator{\Image}{Im}

\DeclareMathOperator{\spn}{span}

\DeclareMathOperator{\Sym}{Sym}
\DeclareMathOperator{\trace}{trace}
\DeclareMathOperator{\Proj}{Proj}
\DeclareMathOperator{\Spec}{Spec}

\DeclareMathOperator{\adj}{ad}

\DeclareMathOperator{\End}{End}
\DeclareMathOperator{\Aut}{Aut}

\DeclareMathOperator{\Lie}{Lie}

\DeclareMathOperator{\Mat}{Mat}
\DeclareMathOperator{\GL}{GL}
\DeclareMathOperator{\gl}{\mathfrak{gl}}
\DeclareMathOperator{\PGL}{PGL}

\DeclareMathOperator{\U}{U}
\DeclareMathOperator{\gu}{\mathfrak{u}}
\newcommand{\uu}{\gu}
\DeclareMathOperator{\SU}{SU}
\DeclareMathOperator{\su}{\mathfrak{su}}
\DeclareMathOperator{\PU}{PU}
\DeclareMathOperator{\Vol}{Vol}

\def\qsing#1/#2(#3){\frac{#1}{#2}(#3)} % to get 1/r(a,b,c,...,d) notation
\def\<#1>{\langle#1\rangle}
\newcommand{\ip}[2]{\langle#1,#2\rangle}

\newcommand{\inv}{{}^{-1}}
\newcommand{\trans}{{}^t\kern -0.2em}

\newcommand{\sst}[1]{#1^{\protect\text{ss}}}
\newcommand{\st}[1]{#1^{\protect\text{s}}}
\newcommand{\aftersub}{\nopagebreak}
\newcommand{\map}[5]{$$\begin{array}{rccc} 
#1\colon & #2 & \longrightarrow & #3 \\ 
         & #4 & \longmapsto & #5
\end{array}$$}
\newcommand{\function}{\map}

\newcommand{\mapeq}[6]{\begin{equation}
\begin{array}{rccc} 
#1\colon & #2 & \longrightarrow     & #3 \\ 
         & #4 & \longmapsto & #5
\end{array}
\protect\label{#6}
\end{equation}}

\newcommand{\corresp}[4]{$$\begin{array}{ccc} 
#1 & \to     & #2 \\ 
#3 & \mapsto & #4
\end{array}$$}

\numberwithin{equation}{section}
\def\*#1:#2*{\S\ref{sec:#1}.\ref{sec:#1:#2}}

        {\begin{list}{}{%
                \settowidth{\labelwidth}{{\emph{#1:}}}%
                \setlength{\leftmargin}{\labelwidth+\labelsep}}}%
        {\end{list}}
\newenvironment{entry}
        {\begin{list}{}%
                {%
                  \setlength{\labelwidth}{12mm}%
                  \setlength{\leftmargin}{14mm}%
                }%
        }%
        {\end{list}}
\newlength{\Mylen}
\newcommand{\Lentrylabel}[1]{%
        \settowidth{\Mylen}{\emph{#1}}%
        \ifthenelse{\lengthtest{\Mylen > \labelwidth}}%
                {\parbox[b]{\labelwidth}%     term > labelwidth 
                        {\makebox[0pt][l]{\emph{#1}}\\}}%
                {\emph{#1}}%                 term < labelwidth 
        \hfil\relax}

\newcommand{\Mentrylabel}[1]%
        {\raisebox{0pt}[1ex][0pt]{\makebox[\labelwidth][l]%
                {\parbox[t]{\labelwidth}{\hspace{0pt}\emph{#1:}}}}}
        {\begin{entry}}%
        {\end{entry}}

\begin{document}
\title[\runningheadstring]{\titlestring\footnote{Maths Subject Classification \msc}}
\author{\name}
\email{sacha@kurims.kyoto-u.ac.jp}
\address{Research Institute for Mathematical Sciences\\ Ky\=oto University\\ 
  Oiwake-ch\protect\=o\\ Kitashirakawa\\ Saky\protect\=o-ku\\ Ky\=oto
  606-01\\ Japan}

\date{2 October 1996} 

\begin{abstract}

\hyphenation{trans-la-tion}

Let~$\Gamma$ be a finite group acting linearly on~$\C^n$, freely
outside the origin, and let $N$ be the number of conjugacy classes of
$\Gamma$ minus one.

A construction of Kronheimer~\cite{kron:ale} of moduli spaces
$X_\zeta$ of translation-invariant $\Gamma$-equivariant instantons on
$\C^2$ is generalised to $\C^n$.

The moduli spaces $X_\zeta$ depend on a parameter~$\zeta\in\Q^N$.  The
following results are proved: for $\zeta=0$, $X_0$ is isomorphic to
$\C^n/\Gamma$; if $\zeta\neq 0$, the natural maps $X_\zeta\to X_0$
are partial resolutions.  The moduli $X_\zeta$ are furthermore shown
to admit K\"ahler metrics which are Asymptotically Locally Euclidean
(ALE).

A description of the singularities of $X_\zeta$ using deformation
complexes is given, and is applied in particular to the case
$\Gamma\subset\SU(3)$.  It is conjectured that for general $\Gamma$
and generic $\zeta$ that the singularities of $X_\zeta$ are at most
quadratic.  When $\Gamma\subset\SU(3)$ a natural holomorphic 3-form is
constructed on the smooth locus of $X_\zeta$, which is conjectured to
be non-vanishing.  The morphims $X_\zeta\to X_0$ are expected to be
crepant resolutions and $X_\zeta$ to be smooth for generic choices of
the parameter $\zeta$.  Related open problems in higher-dimensional
complex geometry are also mentioned.

The paper has a companion paper~\cite{sacha:flows} which identifies
the moduli $X_\zeta$ with representation moduli of McKay quivers, and
describes them completely in the case of abelian groups.

\end{abstract}
\maketitle
\tableofcontents

% ale-main.tex
%

\setcounter{section}{-1}
\section{Introduction}
\label{sec:intro}

This paper is concerned with affine {\em orbifold singularities},
namely with singularities of the type $X=\C^n/\Gamma$ for $\ga$ a
finite group acting linearly on $\C^n$.

More precisely, this paper gives a method for constructing {\em
  partial resolutions\/} of $X$, namely birational morphisms $Y\to X$
which are isomorphisms over the smooth locus of $X$.

\subsection{Background}
\label{sec:intro:back}

The method in question was first introduced by
Kronheimer~\cite{kron:ale}.  It can be described in various ways,
depending on one's point of view.  One description (although maybe not
the most straight-forward) is to construct moduli spaces $X_\zeta$ of
instantons on the trivial bundle $\C^n\times R\to\C^n$. Here $R$
denotes the regular representation space for the group $\Gamma$, and
$\zeta$ is a linearisation of the bundle action.  The instantons are
required to satisfy Hermitian-Yang-Mills-type equations, as well as
additional $\Gamma$-equivariance and translation-invariance
properties.

In Kronheimer's case, $\Gamma\subset\SU(2)$, and the moduli spaces
$X_\zeta$ can in fact be viewed as hyper-\kah\ quotients. Kronheimer
shows that $X_0$ is isomorphic to $\C^2/\Gamma$, that the natural maps
$X_\zeta\to X_0$ are partial resolutions for $\zeta\neq 0$, and that
indeed $X_\zeta$ coincides with the minimal resolution of
$\C^2/\Gamma$ for generic choices of $\zeta$.  Furthermore, $X_\zeta$
inherit natural hyper-\kah\ metrics on their non-singular locus, which
are shown to be Asymptotically Locally Euclidean (ALE): they
asymptotically approximate the Euclidean metric at infinity (up to
terms vanishing with the inverse of the fourth power of the radial
coordinate).

\subsection{Main Results}
\label{sec:intro:main}

In the present case $n$ is any integer greater than or equal to $2$,
and $\Gamma\subset\U(n)$ is assumed to act on $\C^n$ freely outside
the origin for any $n\geq 2$, which means that $X=\C^n/\Gamma$ has an
isolated singularity.\footnote{This is for the purpose of simplicity
  --- the method would seem to be applicable to the general case with
  some modifications} As a result, the moduli $X_\zeta$ are only {\em
  \kah \/} rather than hyper-\kah\ quotients (in actual fact they are
more conveniently described in term of \git).  The main result is
\begin{nonumberthm}[\cf Thms.~\ref{thm:X0free}, \ref{thm:partial_resolutions} and \ref{thm:ale} in the text]
  Let $\Gamma$ act linearly on $\C^n$ and freely outside the origin
  and let $X_\zeta$ be the moduli spaces constructed in
  Section~{\rm \ref{sec:setup:moment}}.
  
  Then $X_0$ is isomorphic to $X=\C^n/\Gamma$, and for $\zeta\neq 0$,
  the natural morphisms $X_\zeta\to X_0$ are partial resolutions.
  Furthermore, the inherited \kah\ metrics on the smooth loci of
  $X_\zeta$ are Asymptotically Locally Euclidean in the sense
  of~\cite{kron:ale}.
\end{nonumberthm}

\subsection{Two Conjectures}
\label{sec:intro:discussion}

The final sections of the paper discuss and develop two conjectures.

\begin{nonumberconj}[\cf Conj. \ref{conj:formal}]
  The singularities of $X_\zeta$ (for generic $\zeta$, say) are at
  most quadratic algebraic.
\end{nonumberconj}

This is a common occurrence for moduli spaces of this
kind~\cite{nadel:quadratic,gold_mill:flat,gold_mill:invariance}.  Its
proof can be reduced to proving the formality of a certain
differential graded Lie algebra (DGLA) by the methods
of~\cite{gold_mill:invariance}.  This is done in
Section~\ref{sec:defcplx}, where the singularities of $X_\zeta$ are
described in terms of deformation complexes~\cite{ahs,don_kron:4mfds}.
The concept of formality is explained, and it is suggested that the
complex relevant to $X_\zeta$ may be formal, in a way similar to
Tian~\cite{tian:smoothness} and Todorov's~\cite{todo:weil-peterson}
work.  This would imply that the singularities of $X_\zeta$ at most
quadratic algebraic.  This conjecture is also checked by computer for
low order abelian groups in~$\U(3)$ using the methods in the companion
paper~\cite{sacha:flows}.

\begin{nonumberconj}[\cf Conj.\ \ref{conj:su3} in the main text]
  If $\Gamma\subset\SU(3)$, the morphisms $X_\zeta\to X_0$ are
  crepant, and if $\zeta$ is generic, $X_\zeta$ is smooth and its
  Euler number is equal to the orbifold Euler number of $X_0$ as
  defined in~\cite{dhvw:i}.
\end{nonumberconj}

The fact that $X_\zeta$ has at most quadratic singularities has been
verified for the abelian subgroups of order less than~11.  The
smoothness of $X_\zeta$ has been verified in the abelian cases $\qsing
1/3(1,1,1)$, $\qsing 1/6(1,2,3)$, $\qsing 1/7(1,2,4)$, $\qsing
1/8(1,2,5)$, $\qsing 1/9(1,2,6)$, $\qsing 1/10(1,2,7)$ and $\qsing
1/11(1,2,8)$.  Both these verifications were done by a brute-force
listing of singularities of $X_\zeta$ for all possible $\zeta$, using
the methods given in the companion paper~\cite{sacha:flows}.

The cases $\Gamma\subset\SU(n)$ present a particular interest.  The
problem of constructing a crepant resolution of $C^3/\Gamma$ with the
same orbifold Euler number was only recently
completed~\cite{mar_ols_per,roan:mirror_cy,mark:res_168,roan:res_a5,ito:trihedral,roan:calabi-yau}.
For the case $\C^4/\Gamma$, one can obtain some interesting analogous
results if one considers terminalisations rather than
resolutions~\cite{sacha:sl4}.
 
In Section~\ref{sec:su3}, a natural holomorphic 3-form is constructed
on the smooth locus of $X_\zeta$: this is conjectured to be
non-vanishing.  Its norm is shown to be constant if and only if the
induced metric on $X_\zeta$ is Ricci-flat (which does not usually turn
out to be the case, however).

\subsection{Related Questions}
\label{sec:intro:related}

This paper has a companion paper~\cite{sacha:flows} in which the
moduli $X_\zeta$ are identified with representation moduli of McKay
quivers.  This allows one to explicitly describe the case of
Abelian~$\Gamma$ in terms of ``flows'' on the McKay quiver.  Explicit
computations are carried out for groups of low order, and the
conjectures about the smoothness and the triviality of the canonical
bundle are checked (by brute-force computer calculations) for abelian
subgroups of $\SU(3)$ of order less than or equal to 11.

Many questions are left open by the present work, besides the
conjectures already mentioned.  For instance, do the birational models
$X_\zeta$ of $\C^n/\Gamma$ possess any special properties with regards
to their singularities (are they terminal?) Are all terminal models
for a given 3-fold singularity obtained by this construction?  What is
the relationship between the different $X_\zeta$? Are they related by
flips/flops? Is it possible, by choosing very special values of
$\zeta$, to produce blowups $X_\zeta\to X_0$ which are interesting
from the point of view of higher dimensional geometry, for instance,
for the construction of flips in dimensions~4 and greater?

\subsection{Methods}
\label{sec:intro:methods}

The methods used include \git, \kah\ quotients, and elementary theory
of the moduli of bundles, the necessary aspects of which are reviewed
in Sections~\ref{sec:git}, \ref{sec:hym} and~\ref{sec:setup}.
Furthermore, the same construction is presented under different angles
with the intention that the reader who is familiar with one of them
(or with Kronheimer's work~\cite{kron:ale}) will be able to follow the
discussion easily.

The later sections devoted to the various conjectures raised by the
main results touch on the theory of deformation complexes, and
concepts of Kuranishi germs, formality, and so on. Some background is
also provided, although not as extensive as to be able to describe it
as ``self-contained'', given the conjectural nature of the material.

\subsection{Outline}
\label{sec:intro:outline}

The outline of this paper is as follows.

Section~\ref{sec:git} reviews material regarding \git\ quotients which
is necessary to define the moduli $X_\zeta$.  No essentially new
material is involved, although the formulation of some of the results
may be un-familiar to non-specialists.

Section~\ref{sec:hym} review material concerning moduli of
Hermitian-Yang-Mills connections. This is not essential to the
understanding of the main results, although some familiarity is
desirable for the understanding of Section~\ref{sec:defcplx}.

Section~\ref{sec:setup} deals with the definition and construction of
the moduli $X_\zeta$.

Section~\ref{sec:partial} gives the proof that $X_0$ is isomorphic to
$X$ and that $X_\zeta\to X_0$ are partial resolutions.

Section~\ref{sec:ale} proves that the induced metrics on $X_\zeta$ are 
ALE.

Section~\ref{sec:defcplx} contains the discussion of the singularities 
of $X_\zeta$ in the language of deformation complexes.  This includes
the conjecture that the singularities of $X_\zeta$ are at most
quadratic algebraic.

Section~\ref{sec:su3} deals with the case $\Gamma\subset\SU(3)$, 
the construction of the holomorphic three-form on the non-singular
locus of $X_\zeta$ and conjecture~\ref{conj:su3}.

\subsection{Acknowledgments}
\label{sec:intro:ack}

The present paper and its companion paper~\cite{sacha:flows} consist
mostly\footnote{Minor portions have been rewritten to include
  references to advances in the field made since then.} of excerpts
of my D.Phil.\ thesis~\cite{sacha:thesis}, and I wish to acknowledge
the University of Oxford and Wolfson College for their hospitality
during its preparation.  I am grateful to the Rhodes Trust for
financial support during my first three years, and to Wolfson College
for a loan in my final year.  The conversion from thesis to article
format was done while I was a Research Assistant in RIMS, Kyoto.

I also take the opportunity to thank my supervisors Peter Kronheimer
and Sir~Michael Atiyah who provided me with constant advice,
encouragement and support and whose mathematical insight has been an
inspiration.  I also wish to thank William Crawley-Boevey, Michel
Brion, Gavin Brown, Jack Evans, Partha Guha, Katrina Hicks, Frances
Kirwan, Alistair Mees, Alvise Munari, Martyn Quick, David Reed, Miles
Reid, Michael Thaddeus, and, last but not least, my parents and
family.

\section{Geometric Invariant Theory}
\label{sec:git}

This section recalls the geometric invariant theory of affine
varieties and proves some results which shall be needed in the sequel.
These results have been included here, because, although they are
well-known to the experts, no elementary treatment
exists.\footnote{For the generalisation of these results to arbitrary
  quasi-projective varieties, see~\cite{thaddeus:git_flips,dolg_hu}.}

\subsection{Linearisations and GIT quotients}
\label{sec:git:lin}

Let~$G$ be a reductive group acting linearly on a complex affine
variety~$X$.  In this situation,\footnote{When $X$ is only
  quasi-projective the definition of a linearisation involves
  specifying an ample bundle over $X$ as well as a lift of the action
  to the bundle.}\/ a ($G$-)\emph{linearisation} is a lifting of the
$G$-action to the trivial line bundle $L\to X$.  Such a linearisation
is determined completely by the action of $G$ on the fibres of $L$,
namely by a character $\zeta\colon G\to \CX$.  For every character
$\zeta$,  denote by $L_\zeta$ the trivial bundle endowed with the
corresponding linearisation.  The space of $G$-invariant sections of
$L_\zeta$ is denoted by~$H^0(L_\zeta)^G$.

The \git\ (GIT) quotient of $X$ by $G$ with respect to $\zeta$ is defined by
$$X\gitquot\zeta G :=\Proj \bigoplus_{k\in\N}H^0(kL_\zeta)^G.$$  

\begin{example}
  \label{example:zeta_zero}
  If $\zeta=0$ is the trivial character, then the corresponding
  quotient $X\gitquot 0G$ coincides with the usual affine GIT quotient
  $X\gitquot{}G$.  In fact, suppose that $X=\Spec R$ for a finitely
  generated ring $R$ and let $z_0$ be a coordinate in the fibre of
 ~$L_0$.  Then $\oplus_{k\in\N}H^0(kL_0)^G=R[z_0]^G=R^G[z_0]$, so
  taking Proj gives: $$X\gitquot 0 G =\Proj R^G[z_0]=\Spec
  R^G=X\gitquot{} G,$$ which is the usual affine GIT quotient.
\end{example}

\subsection{Stability and Extended $G$-equivalence}
\label{sec:git:stab}

The GIT quotient can be obtained by first restricting attention to the
open set $\sst X(\zeta)\subseteq X$ of so called \emph{semi-stable}
points.  A point $x$ in $X$ is called \emph{semi-stable} (with respect
to $\zeta$) if there exists a $G$-invariant section of $kL_\zeta$ (for
some $k$ in $\N$) which is non-vanishing at~$x$.  As a set,
$X\gitquot\zeta G$ is the quotient of $\sst X(\zeta)$ by the
\emph{extended $G$-equivalence} relation induced by the closure of the
$G$-orbits:
$$x\sim y \iff \clos{Gx}\cap\clos{Gy}\neq \emptyset.$$ Thus the
$G$-invariant quotient map $\sst X(\zeta)\to X\gitquot\zeta G$ for the
equivalence relation can map several $G$-orbits to the same point.

For most of the points, this does not happen, however.  This is
because the closure of an open orbit is obtained by adding orbits of
smaller dimension, so since the dimension of the orbit is a lower
semi-continuous function on $X$, it follows that there is an open
subset $\st X(\zeta)\subseteq\sst X(\zeta)$ of points which have
full-dimensional closed $G$-orbits --- the so-called \emph{stable}
points --- and there is a one-one correspondence between the orbits of
$G$ in $\st X(\zeta)$ and their images in the GIT quotient.  In other
words, the GIT quotient contains, as an open set, the geometric quotient
$\st X(\zeta)/G$.

\subsection{Quotients for non-trivial linearisations}
\label{sec:git:gen}

Example~\ref{example:zeta_zero} showed that the GIT quotient for the
trivial linearisation coincides with the affine GIT quotient.  The
following theorem uses the notion of stability to show that the
quotients for non-trivial $\zeta$ are closely related to the affine
quotient.

\begin{thm}
  \label{thm:partial_res}
  The GIT quotients admit projective morphisms
  $$\rho_\zeta\colon X\gitquot\zeta G\to X\gitquot 0G$$
  which are isomorphisms over
  $$\rho_\zeta^{-1}(\st X(0))\to\st X(0).$$
\end{thm}
\begin{proof}
  Let $X=\Spec R$ for some finitely generated ring $R$, and let
  $z_\zeta$ be a complex coordinate in the fibre of~$L_\zeta$.  The
  GIT quotient $X\gitquot\zeta G$ is given by taking Proj of the
  $G$-invariant part of $R[z_\zeta]$, where $R[z_\zeta]$ is to be
  considered as an algebra graded by the powers of~$z_\zeta$.  The
  previous example showed that when $\zeta$ is zero, $\Proj
  R^G[z_0]=\Spec R^G$.
  
  For a general non-zero $\zeta$, $R[z_\zeta]^G\neq R^G[z_\zeta]$.
  However, the {\em degree-zero\/} part of $R[z_\zeta]^G$ is
  always~$R^G$, and this shows~\cite[Example~II.4.8.1 and
  Cor.~II.5.16]{hart:ag} that $X\gitquot\chi G$ is projective over
  $X\gitquot{}G=\Spec R^G$.

  Finally, the map $X\gitquot\zeta G\to X\gitquot 0G$ comes from the
  descent of the composition $\sst X(\zeta)\hookrightarrow \sst X(0)
  \to X\gitquot 0G$ and this is one-one whenever $\sst X(0) \to
  X\gitquot 0G$ is, so the last statement of the theorem follows.
\end{proof}
\begin{rmk}
  \label{rmk:partial_res}
  If $X$ contains a $0$-stable point, then $\st X(0)$ is open and
  non-empty, so dense in~$X$.  Its image $\st X(0)\gitquot 0G$ is open
  and dense in $X\gitquot 0G$, and therefore, $\rho_\zeta\colon
  X\gitquot\zeta G\to X\gitquot 0G$ is an isomorphism on a dense open
  subset, i.e.\ $\rho_\zeta$ is birational.
\end{rmk}

%\section{-}
\section{Moduli of \hym\ Connections}
\label{sec:hym}

This section recalls basic background concerning \hym\ connections and
the construction of their moduli, following~\cite{don_kron:4mfds}.
The construction of $X_\zeta$ will appear as a specialisation of the
material in this section.  However, the reader wishing to go straight
to the point can skip this section and find a self-contained
construction of $X_\zeta$ in section~\ref{sec:setup}.

\subsection{Connections over Symplectic Manifolds}
\label{sec:hym:symplectic}

Suppose $X$ is a compact symplectic manifold $(X,\omega)$ of dimension
$2n$ and let $E$ be a complex vector bundle over~$X$.  The bundle of
infinitesimal automorphisms of $E$ will be denoted $\gl(E)$ or~$\End
E$.

\subsubsection{Connections}

Consider connections on $E$, namely linear maps $\nabla\colon
\Omega^0_{X}(E) \to \Omega^1_{X}(E)$ which satisfy the Leibniz condition.
Any connection on $ E$ can be expressed in a local \nbd\ $U\subset X$
as
\begin{equation}
d_\alpha= d+\alpha,
  \label{eq:dalpha}
\end{equation} for $\alpha\in\Omega^1_U(\End  E)$.  On the other hand, if one considers the difference of two connections, one gets a {\em
  global\/} one-form with values in~$\End E$.

\subsubsection{Hermitian Structure}

Let $\bdh$ be a \hip\ on the fibres of $E$ and denote by $\lieu(
E)\subset\gl( E)$ the real sub-bundle of unitary automorphisms
determined by~$\bdh$.  The connection $\nabla$ is said to be
\emph{compatible} with the Hermitian structure if
$$d \bdh(s,t) =\bdh(\nabla s,t)+\bdh(s,\nabla t).$$ Such a connection
will have local one-form representatives $\alpha\in\Omega_U^1(\lieu(
E))$;  the space $\cA$ of all such connections is an infinite-dimensional
affine space modeled on~$\Omega_X^1(\lieu( E))$.

\begin{rmk}
  Why does one fix a Hermitian structure on $E$?  One reason is
  because fixing a Hermitian structure on $E$ amounts, by Chern-Weil
  theory, to fixing the topological invariants of the connections: for
  any Hermitian connection $d_\alpha$ on $E$, the Chern polynomials
  $c_1(E)$ and $c_2(E)-c_1(E)^2$ are represented respectively by
  $\frac{i}{2\pi}\trace F_\alpha$ and $\trace F_\alpha^2$.
\end{rmk}

\subsubsection{Symplectic Structure and Gauge Group}

The space $\cA$ has a symplectic structure defined by
$$\bomega(a,b):=\int_X\trace(a\wedge b)\wedge\omega^{n-1},$$ for
$a,b\in\Omega_X^1(\lieu( E))$ tangent vectors to~$\cA$.

The \emph{gauge group} $\cG$ is the group of automorphisms of $ E$ which
respect the Hermitian structure in the fibres and cover the identity
map of $X$.  It acts on $\cA$ by $$g\cdot\nabla:= g\nabla
g^{-1}$$ (the condition that $g$ is unitary ensures that the new
connection is compatible with $\bdh$) and preserves the symplectic
form~$\bomega$.

\subsubsection{Moment Map}

The Lie algebra of the gauge group is~$\Lie
\cG=\Omega^0_X(\lieu( E))$.  The moment map for the action of the
gauge group is given by~\cite[Prop.\ 6.5.8]{don_kron:4mfds}
\mapeq\mu\cA{\Omega^0_X(\lieu( E))^*}\alpha{s\mapsto \int_X \trace
  s F_\alpha\wedge\omega^{n-1}}{eq:mucA}
Note that the $\cG$-equivariance follows because the curvature transforms as a tensor under gauge transformations.

\subsection{Connections over \kah\ manifolds, Holomorphic Structures and the \hym\ condition}
\label{sec:hym:kahler}

Now suppose that $X$ is in fact a \kah\ manifold. 

\subsubsection{Complex Structure and the Decomposition of Curvature}

There is a natural decomposition $\Omega^1_{X}(\End E) =
(\Omega^{1,0}_{X}\oplus \Omega^{0,1}_{X})\otimes \Omega^0(\End E)$,
and if a connection is expressed in a local holomorphic frame
$\{z_i\}$ according to~\eqref{eq:dalpha} it takes the form
$$\alpha = \sum_i \alpha_i dz_i - \alpha_i^* d \bar z_i,$$ for
$\alpha_i$ smooth sections of~$\End E$.

The local connection $d_\alpha$ splits into a sum of $(1,0)$ and
$(0,1)$ parts $\partial_\alpha+\dbar_\alpha$ given by:
\begin{align}
  \partial_\alpha &= \partial + \sum_i \alpha_i dz_i\\
 \dbar_\alpha &=  \dbar - \sum_i \alpha_i^* d\zbar_i.
\end{align}
The space $\cA$ becomes a (flat) \kah\ manifold when tangent vectors
are identified with their $(0,1)$ parts, or in other words, when
connections are represented by their $\dbar_\alpha$ operators.  The
holomorphic tangent space to $\cA$ is of course isomorphic to
$\Omega^{0,1}_{X}(\End E)$.

The $(1,1)$ and $(2,0)$ parts of the curvature ${F_\alpha}$ of
$d_\alpha$ are given by
\begin{align}\label{eq:F11}
  F_\alpha^{1,1} &= \sum_{i, j} \left(\frac{\partial\alpha_j}{\partial\zbar_i} - \frac{\partial\alpha_j^*}{\partial z_i}
  -[\alpha_i,\alpha^*_j]\right)dz_i\wedge d\zbar_j,\\ \intertext{and}
\label{eq:F20}
F_\alpha^{2,0} &= \sum_{i,j} \left(-\frac{\partial\alpha_i}{\partial
  {z_j}} +\frac{1}{ 2}[\alpha_i,\alpha_j]\right)dz_i\wedge dz_j,
\end{align}
with $F_\alpha^{0,2}$ equal to minus the Hermitian adjoint of
$F_\alpha^{2,0}$.

\subsubsection{The Hermitian-Yang-Mills condition}

A connection on $E$ is called \emph{\hym} (\emph{HYM}) if the inner
product of its curvature with the \kah\ form $\omega$ is a central
element of $\Omega^1_X(\uu(E))$.  This is in fact a moment map
condition: using the identity
$$F_\alpha\wedge\omega^{n-1}=\frac{1}{ n}\ip{F_\alpha}\omega \omega^n =: \frac{1}{ n}(\Lambda F_\alpha)\omega^n,$$
the map in equation~\eqref{eq:mucA} becomes 
$$\mu(\alpha)(s)=\frac{\Vol(X)}{(n-1)!} \trace(s\Lambda F_\alpha),$$
so, embedding the Lie algebra of $\cG$ its dual in the usual way, we
see that the moment map becomes a constant multiple of 
\map{\mu^*}\cA{\Omega^0_X(\uu(E))}\alpha{\Lambda F_\alpha.}

The moduli space of HYM connections is the \kah\ quotient
$\mu^{*-1}(0)/\cG$.  

\subsubsection{Holomorphic Bundles}

Suppose instead that  the Hermitian connection $d_\alpha$ is required to
induce a holomorphic structure on $E$.  By the Newlander-Nirenberg
theorem, prescribing such a structure is exactly equivalent to
specifying a connection $d_\alpha$ which is
\emph{integrable}, \ie whose $(0,1)$-part is such
that $\dbar_\alpha\circ \dbar_\alpha=0$.  This is equivalent to the
condition that $F_\alpha$ is of type $(1,1)$ and gives a \kah\ 
subvariety ${\cA^{1,1}}\subset\cA$.  This variety parametrises all
the possible holomorphic structures which can be put on the $C^\infty$
bundle $E\to X$.

% \subsubsection{Complexified Gauge Group}

The action of $\cG$ on $\cA$ extends to an action of its complexification
$\cG^\C$, which can be thought of naturally as the group of all
general linear automorphisms of $E$ covering the identity map on $X$.
Put $\tilde g:=(g^*)^{-1}$ and let
\begin{align}
   g\cdot\dbar_\alpha &:= g\dbar_\alpha g^{-1},\\ 
   g\cdot\partial_\alpha &:= \tilde g\partial_\alpha \tilde g^{-1}.
\end{align}

This action of $\cG^\C$ preserves the space $\cA^{1,1}$, and its
orbits are equivalence classes of holomorphic bundles.  To get a nice
moduli space (a quasi-projective variety), one must restrict to
the so-called semi-stable bundles, or in other words, consider a GIT
quotient $\cA^{1,1}\gitquot{}\cG^\C$.  A theorem of Uhlenbeck and
Yau~\cite{uhl_yau:hym} states that the moduli space of stable bundles with
the same topological type as $E\to X$ coincides with the moduli space
of Hermitian-Yang-Mills connections on $E$.  This is an
infinite-dimensional version of the correspondence between symplectic
and algebro-geometric quotients.

\begin{rmk}
  One should really use the quotient of $\cG^\C$ by the scalar
  automorphisms, since they act trivially on $\cA$.  The resulting group then
  has a trivial centre so there is only one linearisation of the
  action; it essentially determined by the degree and rank of $E$.
\end{rmk}

\section{Construction of $X_\zeta$}
\label{sec:setup}

A construction of $X_\zeta$ from scratch will be given in this
section.  

Let $Q$ be an $n$-dimensional complex representation of a finite group
\gat.  Average  over the group elements to get  a \hip\ on $Q$
such that $\ga\subset \U(Q)$. Let ${R}$ be the regular
representation of \gat, i.e.\ 
the free \gat-module which is generated over $\C$ by a basis
$\{e_\gamma | \gamma\in \ga\}$, and on which \gat\ acts via the
morphism $\varphi\colon \ga\to\Aut_\C R$ defined by: $$\gamma\cdot
e_\delta := \varphi(\gamma) e_\delta := e_{\gamma\delta}.$$

\subsection{Invariant HYM Connections}
\label{sec:setup:hym}

\begin{note}
  The reader who has not read Section~\ref{sec:hym} or is not interested
in the ``moduli of bundles'' point of view can jump directly
to~\ref{sec:setup:M}.
\end{note}

The construction of $X_\zeta$ is based on a variation on the
construction in the previous section. It consists, roughly
speaking, in applying the construction to the case where the compact
\kah\ manifold $X$ is replaced by the germ of the singularity
$Q/\Gamma$. 

More precisely, start with $Q^*$, the dual vector-space to $Q$, and
$E=Q^*\times R\to Q^*$ the trivial vector bundle with fibre $R$.
Consider the connections on $E$ which are invariant under all
translations in $Q^*$.  These connections are determined by their
value at one point, and they form a finite-dimensional vectorspace
which can be identified with $M=Q\otimes_\C \End_\C R$ by choosing an
isomorphism $\Omega^1_{Q^*}\cong Q$.  The constructions in the
previous section are now valid, because translation invariance
eliminates any problems one might have with the non-compactness of
base space $Q^*$.  Most aspects are indeed a lot simpler: it suffices to set
all the derivatives of the $\alpha_i$ equal to zero, and to ignore any
integrals over the base space and all the formulas remain valid.  

An unusual feature of these invariant connections is that there are
{\em several\/} moduli spaces: the usual \hym\ condition for a
connection states that the contraction of the curvature with the \kah\ 
form should be a central element of the Lie algebra of the gauge
group.  In case of general connections, there is only one
gauge-invariant momentum level set because the centre of the gauge
group $\cG$ is trivial.  In the case of invariant connections the
gauge group that is relevant is a much smaller group $K^\ga$,
consisting of gauge transformations which are invariant with respect
to all translations and the action of $\ga$.  This group consists of 
unitary endomorphisms of $R$ which commute with the action of $\ga$,
and has a non-trivial centre (consisting of the traceless
\gat-endomorphisms $\zeta\colon R\to R$).  The non-triviality of the
centre means that there are several gauge invariant momentum level
sets, and hence several possible moduli $X_\zeta$.

For the sake of the readers unfamiliar with the material in
section~\ref{sec:hym},  the construction of $X_\zeta$ is given in detail
without making any reference to the bundle construction.

\subsection{The Vector Space $\Mga$}
\label{sec:setup:M}

Let $M=Q\otimes_\C \End_\C R$ and let \gat\ act on $\End_\C R$ by
conjugation via $\ga\stackrel{\varphi}\hookrightarrow \Aut R$:
\begin{equation}\label{eq:gamma_action_endR} 
  \gamma \cdot T := \varphi(\gamma)\,T\, \varphi(\gamma)^{-1} 
\end{equation}
This makes $M$ into a \gat-module. Its
\gat-invariant part is denoted ${M^\ga}$:
\begin{equation}\label{eq:definition_M}
  M^\ga := \left( Q\otimes\End R \right)^\ga.
\end{equation}
The spaces $M^\ga$ and $M$ can be described explicitly by choosing a basis
$\{q_l\}_{l=1}^n$ for $Q$, and defining the \emph{components} of
  $\alpha\in M$ by:
\begin{equation}\label{eq:components_alpha}
  \alpha=\sum_{l=1}^n q_l\otimes \alpha_{l}.
\end{equation}
In this way the elements $\alpha\in M$ are identified with
$n$-tuples of linear maps $\alpha_{i}:R\to R$. The elements of 
$M^\ga$
correspond to those $n$-tuples which satisfy the following
equivariance condition
\begin{equation}\label{eq:alpha_equivariance_condition}
  \sum_l \gamma_{kl} \alpha_{l} = \varphi(\gamma)\alpha_{k}\varphi(\gamma)^{-1},\rlap{$\qquad \forall k,\gamma$,}
\end{equation}
where $\bgamma=(\gamma_{kl})$ is the matrix corresponding to 
the
action of the element $\gamma$ on $Q$ with respect to the basis
$\{q_l\}_{l=1}^n$.

\subsection{Symplectic and \kah\ structure.}
\label{sec:setup:symp}

 Endow $R$ with a fixed \hip\ $\<\ ,\ >$ which makes the standard basis 
$e_\gamma$ 
orthonormal. The inner product on $R$ also defines a real structure 
on the 
space $\End_\C R$ of $\C$-linear endomorphisms of $R$ by the 
Hermitian adjoint operation in the usual way: 
$$\ip{T^*x}y := \ip{x}{Ty}. \qquad x,y\in R. $$
  Define a \hip\ ${\bdh}$ on $M$ by
\mapeq{\bdh}{M\times M}{\C}{(\alpha,\beta)}{\sum_i \trace
  (\alpha_i\beta^*_i)}{eq:dfn_omega} The definition of $\bdh$ is independent of the
basis of $Q$ up to a unitary transformations, and restricts to an
inner product on $M^\ga$.
As usual, $\bdh$ induces two forms on the underlying real 
vector-space to $M$: 
\begin{itemize}
\item a non-degenerate symmetric bilinear form ${\bdg=\Re(\bdh)}$ 
called the \emph{Riemannian metric} associated to $\bdh$ 
\map{\bdg}{M\times M}{\R}{(\alpha,\beta)}{\frac{1}{2}\sum_i \trace 
(\alpha_i\beta^*_i+\beta_i\alpha^*_i)} 

\item a non-degenerate skew-symmetric bilinear form
  ${\omega=\Im(\bdh)}$ called the \emph{K\"ahler} form
  associated to $\bdh$ \map{\omega}{M\times
    M}{\R}{(\alpha,\beta)}{\frac{1}{2\sqrt{-1}} \sum_i
    \trace(\alpha_i\beta^*_i-\beta_i\alpha^*_i)}
\end{itemize}
This gives a Riemannian metric $\bdg$ and a K\"ahler form $\omega$ of 
type
$(1,1)$ on $M$ and $M^\ga$, related as usual by
\begin{equation}
\label{eq:omega_g}
\omega(\alpha, \beta )=\bdg( \alpha , i\beta ).
\end{equation}
This makes $M^\ga, M$ and all their complex subvarieties into \kah\ 
varieties.

The group $\GL (R)$ of automorphisms of $R$ acts on $M$ by conjugation 
on $\End R$:
\begin{equation}
\alpha_i\mapsto g\alpha_i g^{-1},\qquad g\in\GL (R).
  \label{eq:action_glr}
\end{equation}
In fact, the scalars act trivially, and the action descends to an
action of $G:=\PGL(R):=\GL(R)/\GL(1)$.

The subgroup $\glg R$ of endomorphisms which commute with the action
of $\ga$ acts on $M^\ga$ and, in the same way, there is a free action
of $G^\ga=\PGL(R)^\Gamma$ on $M^\ga$. A maximal compact subgroup of $G$
(\resp $G^\ga$) is given by $K=\PU(R)$ (\resp $K^\ga:=
\PU^\Gamma(R)$). The compact group $K$ (\resp $K^\ga$) leaves the
\kah\ structure on $M$ (\resp $M^\ga$) invariant.

\subsection{The Variety  $\protect\cnga$ of Commuting Matrices}
\label{sec:setup:cnga}

 Define the following natural map:
\function{\psi}{Q\otimes\End R}{\Lambda^2 Q\otimes \End R}{\sum_k q_k
  \otimes \alpha_{k}}{\sum_{k,l} q_k\wedge q_l[\alpha_{k},\alpha_{l}]}
This definition is independent of the basis of $Q$ (in terms of
connections, it corresponds to calculating the $(0,2)$ part of the
curvature).  Denote the restriction of $\psi$ to $M^\ga$ by the same
letter.   Define
\begin{equation}\label{eq:definition_cN}
  \cN :=\psi\inv(0)\subset M; 
\end{equation} 
it is a cone (i.e.\ it is invariant under multiplication by non-zero
scalars) which is an intersection of quadrics in $M$ given by the
coordinate functions of $\psi$. In the representation of equation~
\eqref{eq:components_alpha}, its points consist of $n$-tuples of
commuting $r\times r$ matrices:
$$\cN = \{(\alpha_1,\dots,\alpha_n): \alpha_i\in \Mat_r(\C),
[\alpha_i,\alpha_j]=0.\}$$ Its $\ga$-invariant part ${\cN^\ga}$
consists of those commuting matrices satisfying the equivariance
condition~\eqref{eq:alpha_equivariance_condition}.

\subsection{Moment Map}
\label{sec:setup:moment}
\aftersub

Consider the vector space $M$ with the hermitian inner product $h$ and
the action of $K$. The Lie algebra of $K$ is isomorphic to $\su R$,
which consists of traceless skew-Hermitian endomorphisms of $R$. Using
the invariant inner product
$$\ip ab := \trace (a b^*)=-\trace ab,$$ identify $\su R$ with
its dual in the usual way. Then the moment map is for the action of
$K$ is
\mapeq{\mu}M{\su
  R}{\alpha}{\sum_{k}[\alpha^*_{k},\alpha_{k}].}{eq:form_mu} (In the
language of of connections, the map $\mu$ corresponds to contracting
the curvature of the connection with the \kah\ form $\omega=\sum_i
dq_i\wedge d\bar{q_i}\in\Omega_{Q^*}^{1,1}$.)  The moment map for the
action of $K^\ga$ on $M^\ga$ is obtained simply by restriction.

The \kah\ quotients  are defined by
\begin{equation}\label{eq:kaehler_reductions}
  X_\zeta := \frac{\mu\inv(\zeta)\cap \cN^\ga}{K^\ga}, \qquad \zeta\in
  \text{Centre}(\su^\ga R).
\end{equation}
As was remarked in section~\ref{sec:git},
to make this definition rigorous, one needs to make sense of the \kah\ 
structure on $X_\zeta$.  One way to do this is by restricting $\zeta$
to take on integral values.  Then, by the correspondence between \kah\ 
and GIT quotients, one has
\begin{equation}
X_\zeta \cong \cN^\ga\gitquot\zeta G^\ga,
  \label{eq:git_xzeta}
\end{equation}
where $\zeta$ on the right-hand side specifies the linearisation of
the action of $G^\ga$ on the trivial line bundle $\cN^\ga\times\C$.

\begin{rmk}
  In the case $\ga\subset SU(2)$, one has $\Lambda^2 Q\cong R_0$ ---
  the trivial representation.  Identifying $\C^2$ with the
  quaternions, the vector space $M^\ga$ becomes a hyper-K\"ahler
  manifold with 3 distinct complex structures $I,J,K$ and
  corresponding associated K\"ahler forms $\omega_I, \omega_J,
  \omega_K$. The map $\psi$ is then a moment map for the complex
  symplectic form $\omega_\C= \omega_J+i\omega_K$ which is itself
  holomorphic with respect to $\omega_I$ and the quotients $X_\zeta$
  are quotients of $M^\ga$ with respect to the hyper-\kah\ moment map
  given by $(\mu,\psi)$, where the second (complex) variable is set to
  zero.  Kronheimer~\cite{kron:thesis,kron:ale} exploits this fact to
  show that $X_\zeta$ are the minimal resolutions of $\C^2/\Gamma$ for
  generic values of $\zeta$.  Furthermore, by varying the level set of
  $\psi$, he obtains universal deformations.  If $\dim M>2$ however,
  $\psi^{-1}(\zeta)$ is not $G^\Gamma$-invariant for non-zero $\zeta$.
\end{rmk}

\begin{rmk}
  \label{rmk:GQaction}
  In fact, there is a further action on $M$ by $\GL Q$ (acting on $Q$
  on the left).  If $\rho\in\glg Q$, then it is easy to see that this
  preserves $M^\ga$, and indeed $\cN^\ga$. Furthermore, the centre
   $Z(\Gamma)$ of $\Gamma$ is a subgroup of $\glg Q$
  which acts trivially because
  of~\eqref{eq:alpha_equivariance_condition}, so there is an action
  of $G':=\glg Q/Z(\ga)$ which the quotients $X_\zeta $
  inherit. The compact subgroup $K':=\ug
  (Q)/Z(\ga)$ acts in a Hamiltonian fashion, and the moment map for
  this is given by \mapeq{\mu'}{X_\zeta}{(\Lie \U^{\ga}
    Q)^*}{[\alpha]}{b\mapsto \sum_{ij}b_{ij}\trace
    \alpha_i\alpha^*_j.}{eq:form_muQ} In general this does not provide
  one with very much information: for instance, if $Q$ is irreducible,
  $G'=\C^*$ and $\mu'(\alpha)$ is the identity endomorphism of $Q$
  times the sum of the norm squared of the $\alpha_i$'s. In the case
  where $\ga$ is abelian, however, $G'$ contains an algebraic torus of
  dimension $n$ acting freely. The components of $\mu'$ are the norm
  squared of the matrices $\alpha_i$. This is exploited in the
  companion paper~\cite{sacha:flows} to obtain a complete description
  of $X_\zeta$ by the methods of toric geometry.
\end{rmk}

%\section{-}
\section{Variation of Quotients and Partial Resolutions}
\label{sec:partial}

The zero momentum quotient $X_0$ is better understood if viewed as a
two-stage construction: first construct a ``universal quotient"
$\cN_0$ by ignoring the $\ga$-equivariance condition (i.e. perform the
same construction with $\ga$ replaced by the trivial group) and then
obtain $X_0$ as its \gat-invariant part.

\subsection{The Universal Quotient}
\label{sec:partial:uni}
\aftersub

Consider taking symplectic quotients of  $\cN$ with respect to the action of 
$K$.  Since the centre of its Lie algebra is trivial, 
there is only {\em one\/} quotient, with  momentum zero: 
$$\cN_0= \frac{ \cN\cap \mu\inv(0)}{K}=\cN\gitquot{} G.$$

\begin{lemma} 
\label{lemma:N0}
The reduction $\cN_0$ is isomorphic to configuration space of 
$r=|\ga|$ points in $Q$:
$$\cN_0 \cong \Sym^r(Q) := {Q^r/ \Sigma_r},$$
where $\Sigma_r$ denotes the permutation group on $r$ letters acting 
component-wise on the Cartesian product $Q^r$.
\end{lemma}
\begin{proof} The proof simply adapts Kronheimer's \cite[Lemma 
  5.2.1]{kron:thesis}.  It is shown that the $K$-orbits in
  $\cN\cap\mu^{-1}(0)$ can be identified in a one-one way with the
  $\Sigma_r$-orbits in $Q^r$.  Let $\alpha\in\cN\cap \mu\inv(0)$ have
  components $\alpha_i$ with respect to a basis $q_i$. The conditions
  $\psi(\alpha)=\mu(\alpha)=0$ give
\begin{align*}
%\label{eq:psizero}
[\alpha_i,\alpha_j] &=0,\qquad \text{for all }i,j\\
%\label{eq:muzero}
\sum_i [\alpha^*_i,\alpha_i] &=0.
\end{align*}
If one denotes by $A_i$ the operator $\adj(\alpha_i)$, one has,
using the Jacobi identity and the above equations:
\[
\sum_i A^*_i A_i(\alpha^*_j) = \sum_i 
[[\alpha^*_i,\alpha_i],\alpha^*_j]+
[[\alpha^*_j,\alpha^*_i],\alpha_i]= 0.
\] The positivity of $A^*_iA_i$ implies $A^*_iA_i(\alpha^*_j)=0$ for
all $i$ and $j$, and hence $A_j(\alpha^*_j)=[\alpha^*_j,\alpha_j]=0$.
Thus $\cN_0$ is the variety of $n$-tuples of normal commuting
endomorphisms of $R$ modulo simultaneous conjugation. Any such
$n$-tuple can simultaneously diagonalised by conjugation by a unitary
matrix. This means that the orbit of an $n$-tuple is determined by the
eigenvalues of its components; more precisely, there are orthonormal
vectors $v_\gamma\in R$ indexed by the elements $\gamma\in\ga$ and
corresponding eigenvalues $\lambda^i_\gamma\in\C$ such that:
\begin{equation}
\label{eq:alpha_eigenvectors}
\alpha_i(v_\gamma)=\lambda^i_\gamma v_\gamma, \qquad \text{for all } i,\gamma.
\end{equation}
This gives $r$ elements 
$$\lambda_\gamma :=\sum_i \lambda^i_\gamma q_i \in Q,$$ which could be
called the {\em eigenvalues\/} of $\alpha$.  These are defined up to a
permutation, because one can always conjugate by an elementary matrix
which permutes the rows of the $\alpha_i$'s.

In geometrical language, denote by~$\Delta\subset M$ the
subspace of $n$-tuples of matrices which are diagonal with respect to
the standard basis~$e_\gamma$ of~$R$.  The unitary automorphism of $R$
which maps $e_\gamma$ to $v_\gamma$ moves $\alpha$ into $\Delta$.  The
slice $\Delta$ can be identified with $Q^r$ by mapping $\alpha\in
\Delta$ to its $r$ eigenvalues (listed in some specified order).  In
this way $\Delta$ inherits an action of $\Sigma_r$, and the
$\U(R)$-orbit of $\alpha\in\mu^{-1}(0)\cap \cN$ intersects $\Delta$ in
a single $\Sigma_r$-orbit.
\end{proof}

\subsection{The Zero-Momentum Quotient}
\label{sec:partial:zero}
\aftersub

\begin{thm}
  \label{thm:X0free}
  If $\Gamma$ acts freely outside the origin, then $X_0\cong Q/\ga$ as
  varieties.
\end{thm}
\begin{proof}  
  The proof consists in showing that the $K^\ga$-orbits in
  $\cN^\ga\cap\mu^{-1}(0)$ can be identified in a one-one way with the
  $\ga$-orbits in $Q$.  Let $\alpha\in\cN^\ga\cap\mu^{-1}(0)$ and let
  $v_1$ be an eigenvector of $\alpha$ with eigenvalue
  $\lambda_1:=\sum_i\lambda^i_1 q_i\in Q$:
  $$\alpha_i(v_1)=\lambda^i_1 v_1,\quad i=1,\dots,n.$$ By equivariance
  of $\alpha$, 
  $$\alpha_i (R(\gamma)v)=(Q(\gamma)\lambda_1)^i(R(\gamma)v), \text {
    for all } \gamma,\text{ and all }i,$$ so the eigenvectors and
  eigenvalues of $\alpha$ are given by
  $$\lambda_\gamma:=Q(\gamma)\lambda_1\text{ and
    }v_\gamma:=R(\gamma)v_1,\quad\text{for all }\gamma$$ i.e. they lie in orbits of $\ga$.
  Since $\Gamma$ acts freely outside the origin, the eigenvalues are
  either all zero or all distinct and non-zero.  In the latter case, the
  eigenvectors therefore form a basis of $R$.  The unitary 
  automorphism of $R$ defined by  $e_\gamma \mapsto v_\gamma$ commutes with the action of $\ga$, so defines an element of $K^\ga$.
 If $\Delta^\ga\subset M^\ga$ denotes the $n$-tuples of
      endomorphisms of $R$ which are diagonal with respect to the
      standard basis $\{e_\gamma\}$ then the automorphism carries $\alpha$ into $\Delta^\ga$.    The map $\alpha\mapsto \sum_i\lambda^i_1 q_i$
      identifies $\Delta^\ga$ with $Q$ in a manner that is compatible
      with the $\ga$-action on both sides. 
      Furthermore, the $K^\ga$-orbit of $\alpha$ intersects
      $\Delta^\ga$ in precisely one $\Gamma$-orbit.  Thus
      $X_0\cong \Delta^\ga/\Gamma\cong Q/\Gamma$.
\end{proof}

The following lemma will be useful in the section on ALE metrics
(cf.~\cite{kron:thesis}).

\begin{lemma}
\label{lemma:iso_flat}
If $\Gamma$ acts freely outside the origin the map
$$\mu^{-1}(0)\cap\cN^\ga/K^\ga\to \Delta^\ga/\Gamma$$ is an isometry
when $\Delta^\ga$ is given the metric it inherits as a subspace of
$M^\ga$, namely the Euclidean metric.  Furthermore, the bundle
$\mu^{-1}(0)\cap\cN^\ga\to X_0$ is flat.
\end{lemma}
\begin{proof}
  The key point is that the subspace $\Delta^\Gamma$ is everywhere
  orthogonal to the orbits of $K^\ga$: a tangent vector to the orbits
  consists of an $n$-tuple of matrices of the form $[\xi,\alpha_i]$
  for some $\xi\in\su^\ga(R)$, and these matrices are always zero on
  the diagonal, so orthogonal to $\Delta^\Gamma$.  

  This shows that the bundle $\mu^{-1}(0)\cap\cN^\ga \to X_0$ is flat,
  and the definition of the quotient metric on $X_0$ implies that the
  map $X_0=\mu^{-1}(0)\cap\cN^\ga/K^\ga\to\Delta^\ga/\Gamma$ is an
  isometry.
\end{proof}

\subsubsection{Case when $\ga$ does not act freely outside the origin}

Let $\alpha\in\cngamuo$ have an eigenvalue $\lambda\in Q$.
If the stabiliser $\Gamma_\lambda$ of $\lambda$ is trivial, then
$\alpha$ has $r$ distinct eigenvalues, corresponding to the elements
of the orbit $\Gamma\lambda$. This determines the components
$\alpha_i$ completely on the whole of $R$.

On the other hand, if $\lambda$ has a non-trivial stabiliser
$\Gamma_\lambda$, then this determines $\alpha$ on the
sub-representation $W_\lambda := \spn \ga\cdot E_\lambda\subset R$, where
$E_\lambda$ is the eigenspace corresponding to $\lambda$.  In fact,
$E_\lambda$ is a representation of the stabiliser subgroup
$\Gamma_\lambda$ and $W_\lambda$  is simply the  representation 
 of $\Gamma$ induced by $E_\lambda$:
$$W_\lambda=\Ind_{\Gamma_\lambda}^\Gamma E_\lambda.$$

If $\dim E_\lambda<|\Gamma_\lambda|$, then $W_\lambda\neq R$ and 
$\alpha$ restricts to an endomorphism of $W^\perp_\lambda$. Let 
$\lambda'$ be an eigenvalue of the restriction; the equivariance 
condition then determines $\alpha$ on the factor $W_{\lambda'}$. Continuing 
in this way, one obtains a decomposition of $R$: 
$$R=W_\lambda\oplus W_{\lambda'}\oplus\dots.$$ 

From this discussion, one obtains the following description of the quotient $X_0$:
\begin{thm}
  \label{thm:X0}
  There is an inclusion $Q/\ga \hookrightarrow X_0$; this inclusion is
  an isomorphism if and only if $\ga$ acts freely on $Q$ outside the
  origin.
\end{thm}
\begin{proof}
  The first statement follows because, for any orbit $\ga\lambda$ in
  $Q$, one can construct an $n$-tuple $\alpha$ of diagonal matrices
  whose $\lambda$-eigenspace has dimension equal to the stabiliser
  $\Gamma_\lambda$. The orbit of such an $\alpha$ under $K^\ga$ consists of
  commuting matrices with eigenvalue $\gamma\lambda$ with multiplicity
  $|\Gamma_\lambda|$, for all $\gamma\in\Gamma$.

  For the second statement, note that the {\em if\/} direction is
  theorem~\ref{thm:X0free}.  The {\em only if\/} direction follows
  because if $\lambda$ is an eigenvalue of $\alpha$ with non-trivial
  stabiliser and with multiplicity one, one can set $\alpha$ to
  be zero on $W^\perp_\lambda$ (since $0$ is a fixed point of $\ga$,
  the equivariance condition~\eqref{eq:alpha_equivariance_condition}
  does not imply the existence of other eigenvalues). 

  In general, $X_0$ corresponds configurations of $r=|\ga|$ points of
  $Q$ which are unions of orbits of $\ga$, and hence give rise to a
  decomposition of $R$ into induced representations
$$R=\bigoplus_i \Ind_{\Gamma_{\lambda_i}}^\Gamma E_{\lambda_i},$$
where $E_{\lambda_i}$ denote the $\lambda$-eigenspace of an element $\alpha$.
\end{proof}
\begin{rmk}
  When $\Gamma$ doesn't act freely outside the origin, the quotient
  $X_0$ may end up containing all sorts of things.  For instance, for
  the group action $\qsing 1/5(0,1,-1)$, the quotient $X_0$ contains a 
  copy of $\C^3/\Z_5$, a copy of $\C^5$, eight copies of $\C^2$, etc...
\end{rmk}

\subsection{Non-zero Momentum and Partial Resolutions}
\label{sec:partial:non-z}

By theorem~\ref{sec:git}.\ref{thm:partial_res}, in the case where
$\zeta$ is integral, there are projective morphisms $\rho_\zeta\colon
X_\zeta\to X_0$ which are isomorphisms over the set of points which
have finite $K^\ga$-stabilisers.

\begin{prop}
  If $\ga$ acts freely outside the origin, The stabilisers of $K^\ga$
  on $\cN^\ga\cap\mu^{-1}(0)$ are trivial everywhere except at
  $\alpha=0$.
\end{prop}
\begin{proof}
  An automorphism $T$ of $R$ which fixes $\alpha\in
  \cN^\ga\cap\mu^{-1}(0)$ must preserve the (simultaneous) eigenspaces
  of $\alpha$.  If $T$ also commutes with the action of $\ga$, its
  action on an eigenvector $v\in R$ determines its action on the
  linear span of the $\ga$-orbit of $v$.  In the case where $\ga$ acts
  freely and $\alpha$ is non-zero this means that $T$ is only allowed
  to multiply each eigenvector by the same non-zero constant --- and
  this constant must be of modulus one if $T$ is unitary.  Such a $T$
  thus corresponds to the identity element in the quotient group
  $K^\ga=\PU^\ga(R)$.
\end{proof}

Applying the theorem about \kah\ quotients, one gets the following
theorem:

\begin{thm}
\label{thm:partial_resolutions}
  If\/ $\ga$ acts on $Q$ freely outside the origin, and $\zeta$ is
  integral, there are projective morphisms $\rho_\zeta\colon
  X_\zeta\to X_0=Q/\ga$ which are isomorphisms outside the set
  $\rho_\zeta^{-1}(0)$.
\end{thm}

\begin{rmk}
  Even in the case that $\Gamma$ does not act freely outside the
  origin, it is likely that there are still birational maps from  $X_\zeta$ to
  the component of $X_0$ which is isomorphic to $Q/\Gamma$ and which are
  isomorphisms outside the singular set.
\end{rmk}

\section{ALE Metrics}
\label{sec:ale}

The quotients $X_\zeta$ inherit a metric $\bdg_\zeta$ from the metric
$\bdg$ on the ambient space $M^\ga$.  This section shows that these
are ALE metrics.

A metric $\bdg$ on a real $m$-dimensional Riemannian manifold $X$ is
called \emph{asymptotically locally Euclidean}
(\emph{ALE}) if there exists a compact subset $C\subset
X$ whose complement $X\setminus C$ has a finite covering
$\widetilde{X\setminus C}$ which is diffeomorphic to the complement of
a ball in $\R^m$, and such that, in the pulled-back coordinates
$x_1,\dots,x_m$ on $\widetilde{X\setminus C}$, $\bdg$ takes the form
\begin{equation}
\bdg_{ij}=\delta_{ij} + a_{ij},
  \label{eq:ale}
\end{equation}
where $|\partial^p a_{ij}|=O(r^{-4-p})$ for $p\geq 0$, where
$r=\sqrt{\sum_i x_i^2}$ denotes the radial distance in $\R^m$ and
$\partial$ denotes the differentiation with respect to the coordinates
$x_1,\dots,x_m$.
\begin{thm}
  \label{thm:ale}
  The metrics on $X_\zeta$ are ALE: for any $\zeta$, there is an
  expansion in powers of $r$
  \begin{equation}
    \bdg_\zeta = \delta + \sum_{k\geq 2} h_k(\theta)r^{-2k},
    \label{eq:expansion}
  \end{equation}
  where $(r,\theta)$ denote polar coordinates in $\R^{2n}\cong Q$. This
  expansion is analytic and may  be differentiated term by term.
\end{thm}
\begin{proof}
  Kronheimer's proof~\cite[Prop.5.5.1]{kron:thesis} goes through with
  the appropriate modifications.  The metric $\bdg_\zeta$ restricted
  to the unit ball $r=1$ is an analytic function of $\zeta$, so admits
  an expansion $$\bdg_{\zeta| r=1} = \sum_\nu f_\nu \zeta^\nu$$ where
  $\nu$ are multi-indices in the coordinates of $\zeta$.  The moment
  map being quadratic homogeneous implies that
  $$\bdg_\zeta(r,\theta)=\bdg_{r^{-2}\zeta}(1,\theta).$$ Hence the
  expansion for $\bdg_\zeta$ takes the form $$\sum_{k\geq 0}
  h_k(\theta)r^{-2k},$$ where the $h_k=\sum_{|\nu|=k}f_\nu\zeta^\nu$
  are analytic functions of the radial coordinates.

  It remains to show that $h_0=\delta$ and that $h_1=0$.  The first
  statement is equivalent to showing that the identification $X_0\to
  Q/\Gamma$ is an isometry. This was done in
  Lemma~\ref{lemma:iso_flat}.

  For the second statement, one must show that the variation of
  $\bdg_\zeta$ with $\zeta$ is zero at $\zeta=0$ in every direction
  $\lambda\in ((\Lie K^\ga)^*)^{K^\ga}$. The metric $\bdg_\zeta$ is determined
  entirely by the \kah\ form $\omega_\zeta$ and the induced complex
  structure $J_\zeta$.  Since the latter is the same for all $\zeta$,
  it is sufficient to prove that $$\partial_\lambda
  \omega_{\zeta|\zeta=0}=0$$ for all $\lambda$.  A general formula for
  the variation of the induced symplectic form is given by Duistermaat
  and Heckman in~\cite{dui_hec:variation}.  Away from the
  singularities, the projection $\mu^{-1}(\zeta)\cap\cN^\ga\to
  X_\zeta$ is a principal $K^\ga$-bundle whose connection is given by
  the Levi-Civita connection for the induced metric on $X_\zeta$. If
  $\Omega_\zeta$ denotes the curvature, regarded locally as an element
  of $\Omega^2_{X_\zeta}\otimes\su^\ga(R)$, then the formula for the
  variation of $\omega_\zeta$ is given by $$\partial_\lambda
  \omega_\zeta = \<\lambda,\Omega_\zeta>.$$ In the present case,
  lemma~\ref{lemma:iso_flat} tells us that $\Omega_0=0$, so the
  variation is zero for $\zeta=0$, and this concludes the proof.
\end{proof}

\section{Deformation Complexes}
\label{sec:defcplx}

The question of the local geometry of the moduli spaces $X_\zeta$ can
be studied using the tools of deformation complexes.

\subsection{Differential Forms and Graded Lie Algebras}
\label{sec:defcplx:diff}

Define the vectorspaces $$M^{p,q} := \Omega^{p,q}_{Q^*}\otimes \End R,$$ for
$p,q \in\N$, whose typical element $\beta$ is of the form
$$\beta = \beta^I_{\bar J} dq_I\wedge d\bar q^{J},$$ where the
summation convention is used for the multi-indices $I,J$ and where, as
usual, $dq_I= \wedge_{i\in I} dq_i$, and $d\bar q^J=\wedge_{j\in J}
d\bar q^j$.  Define the \emph{degree} of $\beta$ to be $\deg
\beta:=|I|+|J|$ and write $(-1)^\beta$ for $(-1)^{\deg\beta}$.  The product of two elements $\alpha,\beta$ is
defined to be
$$\alpha\beta:=\alpha^I_{\bar J}\beta^{I'}_{\bar J'}dq_I\wedge d\qbar^{J}\wedge
dq_{I'}\wedge d\qbar^{J'},
$$
and the bracket of
any two elements $\alpha\in M^{p,q}$ and $\beta\in M^{p',q'}$ is
defined by
\begin{equation}
[\alpha,\beta] := [\alpha^I_{\bar J},\beta^{I'}_{\bar J'}]dq_I\wedge d\bar q^{J}\wedge dq_{I'}\wedge d\bar q^{J'} = \alpha\beta-(-1)^{\alpha\beta}\beta\alpha.
  \label{eq:bracket}
\end{equation}
In the equation above and elsewhere, $(-1)^{\alpha\beta}$ means
$(-1)^{\deg\alpha\deg\beta}$ and not
$(-1)^{\deg\alpha}(-1)^{\deg\beta}$. Writing
$M^r:=\sum_{p+q=r}M^{p,q}$, the algebra $M^*$ inherits the structure of a
graded Lie algebra, namely a graded algebra with a bracket satisfying
$$[M^r,M^s]\subset M^{r+s},$$
(graded) skew-commutativity
$$[\alpha,\beta]=-(-1)^{\alpha\beta}[\beta,\alpha],$$ and the (graded)
Jacobi identity:
$$(-1)^{\alpha\gamma}[\alpha,[\beta,\gamma]]+(-1)^{\beta\alpha}[\beta,[\gamma,\alpha]]+(-1)^{\gamma\beta}[\gamma,[\alpha,\beta]]=0.$$
There are two sub-algebras $M^{*,0}$ and $M^{0,*}$ of $M^*$.  Defining
the adjoint of $\alpha=\alpha^I_{\bar J}dq_I\wedge d\qbar^{J}$ to be
$$\alpha^* := (\alpha^I_{\bar J})^*d\qbar^I\wedge dq_{J},$$
then
$(\alpha\beta)^* = (-1)^{\alpha\beta}\beta^*\alpha^*$ and
$[\alpha,\beta]^*=(-1)^{\alpha\beta}[\beta^*,\alpha^*]$.

The Jacobi identity implies that if $\alpha$ or $\beta$ has odd degree then 
\begin{equation}
[\alpha,[\alpha,\beta]]=\frac{1}{2}[[\alpha,\alpha],\beta].
  \label{eq:aab}
\end{equation}
If $\alpha=\alpha_i d\qbar^i \in M^{0,1}$, and one 
defines
\map{\dbar_\alpha}{M^{p,q}}{M^{p,q+1}}{\beta}{[\alpha,\beta],} then,
using~\eqref{eq:aab}, one sees that the sequences
\begin{equation} \label{eq:Mpx_complex}
   M^{p,*}: M^{p,0}\xrightarrow{\dbar^0_\alpha}M^{p,1}\xrightarrow{\dbar^1_\alpha}M^{p,2}\to\dots
\end{equation}
of vectorspaces are complexes precisely when $[\alpha,\alpha]=0$,
i.e.\ when $\alpha\in \cN$.  Write $H^{p,q}_\alpha$ for the cohomology
groups $H^q(M^{p,*},\dbar_\alpha)$.  If one introduces a metric on
$M^{p,q}$ by using the standard inner product on $\End R$ and making
$$\frac{1}{\sqrt{2^{p+q}}}\{dq_I\wedge d\qbar^{J}\}_{\substack{|I|=p,\\
|J|=q}}$$ 
orthonormal~\cite[p.80]{gri_har:ag}, one can define the adjoint
operator $\dbar^*_\alpha$, and the Laplacian
$\bar\Box_\alpha:=\dbar_\alpha\dbar^*_\alpha+\dbar^*_\alpha\dbar_\alpha$.
Their kernels give harmonic representatives for the cohomology groups
in the usual way $$\cH^{p,q}_\alpha:= \ker \bar\Box^{p,q}_\alpha
\subset M^{p,q}.$$
If $\Lambda\colon M^{p,q}\to M^{p-1,q-1}$ denotes
the operation of contraction with the \kah\ form $\omega= dq_i\wedge
d\qbar^i$, then the definition of the adjoint and the invariance of
the trace under cyclic permutations give
\begin{equation}
\dbar^*_\alpha \beta= -\Lambda[\alpha^*,\beta],
  \label{eq:dstar}
\end{equation}
or, in coordinates,
\begin{equation}(\dbar^*_\alpha \beta)^I_{\bar J} = [\alpha^*_j,\beta^I_{j\bar J}].\label{eq:dstar_coord}
\end{equation}
Writing $\kappa:= dq_1\wedge \dots \wedge dq_n$,
the $n$-th power of $\omega$ is
\begin{equation}
\omega^n=n!(-1)^{(n-1)(n-2)/2}\frac{i^n}{2^n}\kappa\wedge\kappa^*.
  \label{eq:omegan}
\end{equation}

\subsection{Local Description of $X_\zeta$}
\label{sec:defcplx:local}

Under the identification $M=M^{0,1}$, the derivative of the
action~\eqref{eq:action_glr} of $\GL(R)$ on $M$ is given by
$$\partial_\alpha^0\colon M^{0,0}\to M^{0,1}.$$  On the other hand, the
derivative of $\alpha\mapsto F^{0,2}_\alpha$ is (twice) 
$$\dbar^1_\alpha \colon M^{0,1}\to M^{0,2},$$ so the Zariski tangent
space to $\cN\gitquot{}\GL(R)$ at an element $\alpha$ is given by the first cohomology group of the Atiyah-Hitchin-Singer~\cite{ahs} deformation complex
$$ M^{0,0}\xrightarrow{\dbar^0_\alpha}M^{0,1}\xrightarrow{\dbar^1_\alpha}M^{0,2},$$ i.e.\ 
by $H^{0,1}_\alpha$.  From the point of view of the \kah\ quotient,
one can see this as follows: the Zariski tangent space to $X_\zeta$ at
$[\alpha]$ is given by $\ker d\mu(\alpha)\cap\ker d\psi(\alpha)$.  By
definition, the derivative of $\mu$ is dual to the action of $K^\ga$,
so $d\mu(\alpha)=-2\dbar^*_\alpha= 2\Lambda[\alpha^*,\ ]$, as can be
verified by remarking that $\mu(\alpha)=\Lambda[\alpha^*,\alpha]$.
Hence  $\ker d\mu(\alpha)=\Image \dbar^*_\alpha$, so the tangent
spaces indeed coincide.

A local model for $\cN\gitquot{}\GL(R)$ in a neighbourhood of a point
$[\alpha]$ where $\alpha\in\cN$ has trivial stabiliser is given by
solving the equation $F^{0,2}_{a+\beta}=0$ for $\beta$ in the slice
$$\{\beta\in M^{0,1}: \dbar_\alpha^*\beta=0,\|\beta\|\text{
  small}\}.$$ This comes down to solving the system of equations
\begin{align} 
&\dbar^*_\alpha \beta=0\\ 
&\dbar_\alpha \beta + \frac{1}{2}[\beta,\beta]=0,
\end{align}
in a neighbourhood of the origin.  Kuranishi's argument~\cite{kura:complex,kura:newproof} shows that the
solution set is given by the zero set of a map $\Phi\colon
\cH^{0,1}_\alpha\to \cH^{0,2}_\alpha$ whose two-jet at the origin is
given by
\map{\Phi_{(2)}}{\cH^{0,1}_\alpha}{\cH^{0,2}_\alpha}{\beta}{\cH_\alpha([\beta,\beta]).}
where
$$\cH_\alpha \colon M^{0,2} \to \cH^{0,2}_\alpha$$ denotes the
orthogonal projection to the harmonic subspace. Similar statements
hold for $X_\zeta$ and $M^{0,*,\ga}$.

\subsection{Kuranishi Germs and Formality}
\label{sec:defcplx:kura}

This whole discussion can be phrased in more abstract language of
deformation functors and differential graded Lie algebras. Additional
details and background can be found
in~\cite{gold_mill:invariance,gold_mill:fundamental} and~\cite{dgms}.
The algebra $M^{0,*}$ is actually a \emph{differential graded Lie
  algebra} (DGLA) when
endowed with the differential $\dbar_\alpha$.  The metric on $M^{p,q}$
makes it into an \emph{analytic DGLA}, namely a DGLA which possesses
a norm compatible with its differential and bracket, and which induces
what is essentially a Hodge decomposition of its graded pieces with
finite dimensional topological summands $\cH^i$ which are the
analogues of the harmonic forms.  When $\alpha$ has trivial
stabiliser, $X_\zeta$ is locally analytically isomorphic in the
neighbourhood of $[\alpha]$ to the \emph{Kuranishi germ} ${\bdK_M}$
associated to $M^{0,*,\ga}$.   The results so far are stated in the
following theorem.

\begin{thm}
  \label{thm:formal}
  The sequence of vector spaces $(M^{0,*},\dbar_\alpha)$ is a complex
  (and therefore a differential graded Lie algebra) if and only if
  $\alpha\in\cN$.  Furthermore, if $\alpha\in\cN^\ga$, the Zariski
  tangent space to $X_\zeta$ at $[\alpha]$ is isomorphic to the first
  cohomology group of its $\ga$-invariant part
  $$H^{0,1,\ga}_\alpha:=H^1(M^{0,*,\Gamma},\dbar_\alpha)$$ and if
  $\alpha$ has trivial $K^\ga$-stabilisers, then $X_\zeta$ is locally
  isomorphic to its Kuranishi germ
  $$\bdK_{M^{0,*,\ga}}=\{\beta\in M^{0,1,\ga} |
  \dbar^*_\alpha\beta=\dbar_\alpha \beta + \frac{1}{2}[\beta,\beta]=0\}.$$ 
\end{thm}

In general the Kuranishi germ $\bdK_L$ of an analytic DGLA $(L,d)$ is
(Banach analytically) isomorphic to the germ at $0$ of
$$\{\beta\in (\Image d)^\perp\subset L^1 | d\beta+\frac{1}{2}[\beta,\beta]=0\},$$ where $(\Image d)^\perp$ is a fixed complement
of the image of $d$ in $L^1$.  Goldman
and~Millson~\cite{gold_mill:invariance} prove that $\bdK_L$ is
invariant under \emph{quasi-isomorphisms}, namely chains of
homomorphisms of DGLAs\footnote{An important point is that the intermediate
  $L',L'',\dots$ need {\em not\/} have any analytic structure and that the
  intermediate arrows need not preserve any splittings.}
$$L\to L' \leftarrow L'' \to \dots \leftarrow L'''$$ which induce
isomorphisms in cohomology.  When $L$ is quasi-isomorphic to its
cohomology (which is a DGLA when endowed with the zero differential),
$L$ is called \emph{formal} and it follows that $\bdK_L$ is analytically
isomorphic to the quadratic cone $$\bdQ_L:=\{\beta\in\cH^1 :
[\beta,\beta]=0\}.$$

One  way in which this can happen is if
the bracket of two harmonic elements of degree $1$ is harmonic.  This
is the case, for instance for the moduli space of flat
Hermitian-Yang-Mills connections over a compact \kah\ 
manifold~\cite{nadel:quadratic,gold_mill:flat,gold_mill:invariance}.
If in addition, the cup-product on $\cH^1$ is zero, then $\bdK_L\cong
\cH^1$ and the deformation space is a smooth manifold (even if
$\cH^2\neq 0$).  This is the case, for instance for the moduli space
of complex structures over a \emph{Calabi-Yau $n$-fold}, namely a
compact \kah\ manifold with a nowhere vanishing holomorphic
$(n,0)$-form~\cite{gold_mill:invariance}.  These moduli were studied by F.~Bogomolov.  The key fact which implies
the formality of the DGLA and the vanishing of the cup-product in this
case was proved by Tian~\cite{tian:smoothness} and
Todorov~\cite{todo:weil-peterson}.

In the case of the algebra $M^{0,*}$, formula~\eqref{eq:aab} with
$\alpha$ and $\beta$ interchanged shows that the bracket of two
harmonic elements in $\cH^{0,1}_\alpha$ is $\dbar_\alpha$-closed.
However, it does not follow that $\dbar^*_\alpha([\beta,\beta])=0$;
indeed this is easily seen to be false, since
$[\beta,\beta]=2\beta\beta$.  Nevertheless, it does not seem
unreasonable to expect that $M^{0,*,\ga}$ can also be proved to be
formal for generic $\zeta$, maybe by imitating Tian and Todorov's
method.
\begin{conj}
\label{conj:formal}
The differential graded Lie algebra $(M^{0,*,\ga},\dbar_\alpha)$ is formal
for all $\alpha\in\cN^\ga\cap\mu^{-1}(\zeta)$ and generic $\zeta$, and therefore $X_\zeta$ has, for these $\zeta$, at worst
quadratic algebraic singularities.
\end{conj} 

Another conjecture is the following:
\begin{conj}
\label{conj:su3}
    If $\Gamma\subset\SU(3)$, can one imitate the Tian-Todorov proof and
  show that the Kuranishi germ of $(M^{0,*,\ga},\dbar_\alpha)$ is
  isomorphic to $\cH^{0,1,\ga}_\alpha$ for generic $\zeta$, i.e.~that
  $X_\zeta$ is smooth?
\end{conj}

The fact that $X_\zeta$ has at most quadratic singularities has been
verified for the abelian subgroups of order less than~11.  The
smoothness of $X_\zeta$ has been verified in the abelian cases $\qsing
1/3(1,1,1)$, $\qsing 1/6(1,2,3)$, $\qsing 1/7(1,2,4)$, $\qsing
1/8(1,2,5)$, $\qsing 1/9(1,2,6)$, $\qsing 1/10(1,2,7)$ and $\qsing
1/11(1,2,8)$.  Both these verifications were done by exhaustive
listing of singularities of $X_\zeta$ for all possible $\zeta$, using
the methods given in the companion paper~\cite{sacha:flows}.

A different approach is available in the specific case of $\SU(3)$;
this is presented next.

\section{Subgroups of $\protect\SU(3)$  and Cubic Forms}
\label{sec:su3}

Suppose that $\Gamma\subset\SU(3)$.  If $\alpha\in\mu^{-1}(\zeta)$ and
$\beta,\delta\in\cH^{0,1,\ga}_\alpha$, then, as  remarked in the
previous section,
$$\dbar_\alpha[\beta,\delta] =0,$$ but $[\beta,\delta]$ is not in
$\cH^{0,2,\ga}_\alpha$.  However, considerations of type show that it differs
from its harmonic projection by a term $\dbar_\alpha\epsilon$, for some
$\epsilon\in M^{0,1,\ga}$. For $\eta\in\cH^{0,1,\ga}_\alpha$
  \begin{align}
\trace(\eta[\beta,\delta])- \trace(\eta\cH_\alpha([\beta,\delta])) &= \trace(\eta[\alpha,\epsilon]) \notag\\
&= \trace(\epsilon[\eta,\alpha])\notag\\
&= 0,\qquad \text{since } \eta\in\cH^{0,1,\ga}_\alpha.
  \label{eq:trace_harmonic}
\end{align}
This shows that the tensor
\corresp{H^{0,1,\Gamma}_\alpha\otimes H^{0,1,\Gamma}_\alpha\otimes
  H^{0,1,\Gamma}_\alpha}{\C}{(\eta,\beta,\delta)\phantom{{}^{,1,\Gamma}}}{\kappa\trace(\eta\cH_\alpha([\beta,\delta])),}
is totally symmetric on $H^{0,1,\Gamma}_\alpha$ (the
isomorphism $\Omega^{3,3}_{Q^*}\cong\C$ has been used).  An easy polarisation
argument shows that it is completely determined by the corresponding cubic form
\map{\bdC}{H^{0,1,\Gamma}_\alpha}\C\beta{\kappa\trace(\beta([\beta,\beta])).}

\begin{prop}
  \label{prop:cubic}
  The singularity of $X_\zeta$ has no quadratic part if and only if
  $\bdC(\beta)=0$ for all $\beta\in\cH^{0,1,\ga}_\alpha$ and all
  $\alpha\in\mu^{-1}(\zeta)$.
\end{prop}
\begin{proof}
  Suppose $X_\zeta$ has no quadratic part at $[\alpha]$.  Then
  $\Phi_{(2)}(\beta)=\cH_\alpha([\beta,\beta])=0$.  But this implies
  that $\bdC(\beta)=0$ by equation~\eqref{eq:trace_harmonic}.

  Conversely, if $\bdC(\beta)=0$ for all $\beta\in\cH^{0,1,\ga}$ then
  the corresponding totally symmetric tensor vanishes on all triples
  $(\eta,\beta,\beta)$ for all $\eta,\beta\in \cH^{0,1,\ga}_\alpha$.
  Since this is true for all $\eta$, it must be that
  $\cH_\alpha([\beta,\beta])\in\Image\dbar_\alpha$, i.e.\
  $\Phi_{(2)}(\beta)=0$ in $\cH^{0,2,\ga}_\alpha$.
\end{proof}
There is a natural $3$-vector ${\Omega}$ whose value on three elements
of $H^{0,1,\ga}_\alpha$ of is given by
\begin{equation}
  \Omega(\eta,\beta,\delta):=\kappa\trace( \eta\beta\delta).
  \label{eq:Omega}
\end{equation}
This is symmetric under cyclic permutations of the entries, so
decomposes into a totally skew-symmetric part $\Omega_{\text{skew}}$
and a totally symmetric part, which is nothing but the totally
symmetric tensor corresponding to $\bdC$.  The proposition above
implies the
\begin{cor}
  \label{cor:3form}
  If $X_\zeta$ is smooth, then $\Omega$ defines an element of
  $\Omega^{3,0}(X_\zeta)$.
\end{cor}
\begin{conj} 
  The canonical sheaf $\cO_{X_\zeta}(K_{X_\zeta})$ is locally free,
  and is generated by the non-vanishing $(3,0)$-form $\Omega$ when
  $X_\zeta$ is smooth.
\end{conj}
Taking the  wedge of $\Omega$ with its complex conjugate gives
\begin{align}
\label{eq:owo1}
  \Omega\wedge\Omega^*(\eta,\beta,\delta,\eta^*,\beta^*,\delta^*) &=
(\epsilon_{ijk}\trace
\eta_i\beta_j\delta_k)\overline{(\epsilon_{\ibar\jbar\kbar}\trace
\eta_\ibar\beta_\jbar\delta_\kbar)}\kappa\wedge\kappa^*,\\
\label{eq:owo2}
&=\left|\sum^\circ_{ijk}\trace \eta_i\beta_j\delta_k\right|^2\kappa\wedge\kappa^*,
\end{align}
where $\sum^\circ_{ijk}$ denotes the sum over distinct $i,j$ and $k$.
On the other hand, the symplectic form $\omega_\zeta$ on $X_\zeta$ is
simply the restriction of the symplectic form $\omega$ defined
in~\eqref{eq:dfn_omega}, and equation~\eqref{eq:omegan} gives
\begin{equation}
  \omega_\zeta\wedge\omega_\zeta\wedge\omega_\zeta  = \frac{3i}{4}\kappa\wedge\kappa^*.
\end{equation}
Suppose that $(\eta,\beta,\delta)$ are an orthonormal triple in
$T^{1,0}_\alpha X_\zeta$.  Then the value of the coefficient of
$\kappa\wedge\kappa^*$ in~\eqref{eq:owo2} is equal to $\|\Omega\|^2
\|\kappa\wedge\kappa^*\|^{-2}$.   Hence $X_\zeta$ has trivial canonical
bundle if this coefficient is never zero for all
$\alpha\in\mu^{-1}(\zeta)$.  
\begin{lemma}
  \label{lemma:bochner} % nothing to do with him!
  The \kah\ manifold $X_\zeta$ is Ricci-flat if
  and only if the coefficient of $\kappa\wedge \kappa^*$
  in~\eqref{eq:owo2} is constant for all orthonormal triples
  $(\eta,\beta,\delta)$ in $H^{0,1,\ga}_\alpha$ and all
  $\alpha\in\mu^{-1}(\zeta)\cap\cN^\ga$.
\end{lemma}
\begin{proof}
  This follows because if $X_\zeta$ is Ricci-flat, there exists a
  holomorphic $(3,0)$-form $\Omega'$ which is covariant constant on
  $X_\zeta$. Hence $\Omega$ will differ from $\Omega'$ by a
  holomorphic function $f$. Now Liouville's theorem implies that $f$
  is either constant or unbounded. Since $\Omega$ is clearly bounded
  (by $6\kappa\wedge\kappa^*$), $f$ must be constant.
\end{proof}

\subsection{Example}
\label{sec:su3:ex}

Let us work out a specific example.  Consider the group $\Gamma=\mu_3$
of order $3$ acting on $\C^3$ with weights $(1,1,1)$.  The following
configuration of matrices is easily seen to define a point of
$\mu^{-1}(\zeta)\cap\cN^\ga$, where $\zeta=(-|A|^2,|A|^2-|B|^2,|B|^2)$, $(A,B\in\R)$:
\begin{equation}
\label{eq:point1}
  \alpha_1=\begin{pmatrix}
    0 &A &0 \\
    0 &0 &B \\
    0 &0 &0
  \end{pmatrix}d\qbar^1, \quad \alpha_2=\alpha_3=0.
\end{equation}
The tangent space is three-dimensional and is generated by the
following orthonormal elements (recall that $\|d\qbar^i\|^2=2$)
\begin{equation}
  \beta_1=\frac{1}{\sqrt 2}\begin{pmatrix}
    0 &0 &0 \\
    0 &0 &0 \\
    1 &0 &0
  \end{pmatrix}d\qbar^1, \quad
 \beta_i= \frac{1}{\sqrt{2(A^2+B^2)}}\begin{pmatrix}
    0 &A &0 \\
    0 &0 &B \\
    0 &0 &0
  \end{pmatrix}d\qbar^i,
\end{equation}
for $i=2,3$, so this defines a smooth point of $X_\zeta$.  The value
of $\|\Omega\|^2$ at this point is
\begin{equation}
\left|\frac{1.A.B+1.B.A}{2\sqrt{2}(A^2+B^2)}\right|^2\kappa\wedge\kappa^* =
\frac{1}{2}\left(\frac{AB}{A^2+B^2}\right)^2 \kappa\wedge\kappa^*,
  \label{eq:norm1}
\end{equation} and so this is non-zero away from $AB=0$ (which correspond to
non-generic values of $\zeta$).

At the point
\begin{equation}\label{eq:point2}
\alpha_1=\begin{pmatrix}
    0 &A+C &0 \\
    0 &0 &B+C \\
    C &0 &0
  \end{pmatrix}d\qbar^1, \quad
   \alpha_2= \alpha_3=0,
\end{equation}
in $\mu^{-1}(\zeta)\cap \cN^\ga$, the tangent space is still
three-dimensional, with orthonormal generators
\begin{equation}
  \beta_i = \frac{1}{\sqrt{6}}\begin{pmatrix}
    0 &1 &0 \\
    0 &0 &1 \\
    1 &0 &0
  \end{pmatrix}d\qbar^i, \quad i=1,2,3,
\end{equation}
The value of $\|\Omega\|^2$ however is now
\begin{equation}
\left|6\left(\frac{1}{\sqrt 6}\right)^3\right|^2 \kappa\wedge\kappa^*=\frac{1}{6}\kappa\wedge\kappa^*.
  \label{eq:norm2}
\end{equation}

In fact, all the points of $\mu^{-1}(\zeta)\cap\cN^\ga$ are of the
form~ \eqref{eq:point1} or ~\eqref{eq:point2} (modulo permutations of
the indices $1,2,3$)\footnote{See~\cite{sacha:flows} for an
  explanation of why this is so.}.  Thus it has been shown, in a
rather laborious way, that away from certain degenerate values of
$\zeta$, $\Omega$ is non-vanishing on $X_\zeta$ and $K_{X_\zeta}$ is
therefore trivial. In fact, $X_\zeta=\cO_{\PP^2}(-3)$.

Since the coefficient of $\kappa\wedge\kappa^*$ in~\eqref{eq:norm1} is
always smaller than ${1/ 8}$, one also deduces that
$\Omega\wedge\Omega^*$ is not a constant multiple of
$\omega_\zeta\wedge\omega_\zeta\wedge\omega_\zeta$ on any of the
quotients $X_\zeta$, and therefore by lemma~\ref{lemma:bochner} that
the induced metric is never Ricci-flat.

\begin{rmk} The space $\cO_{\PP^2}(-3)$ does have a standard
Ricci-flat metric, as was first noted by Calabi~\cite{calabi}.
\end{rmk}

\providecommand{\bysame}{\leavevmode\hbox to3em{\hrulefill}\thinspace}


\begin{thebibliography}{DHVW85}

\bibitem[AHS77]{ahs}
M.~F. Atiyah, N.J. Hitchin, and I.M. Singer, \emph{Deformations of instantons},
  Proc. Nat. Acad. Sci. U.S.A. \textbf{74} (1977), no.~7, 2662--2663.

\bibitem[Cal79]{calabi}
Eugenio Calabi, \emph{M\'etriques {K\"a}hl\'eriennes et fibr\'es holomorphes},
  Ann. Ec. Norm. Sup. ($4^{\rm e}$ s\'erie) \textbf{12} (1979), 269--294.

\bibitem[DGMS75]{dgms}
P.~Deligne, J.~Griffiths, J.~Morgan, and D.~Sullivan, \emph{Real homotopy
  theory of {K\"a}hler manifolds}, Inventiones Math. \textbf{29} (1975),
  245--274.

\bibitem[DH82]{dui_hec:variation}
J.~J. Duistermaat and G.~J. Heckman, \emph{On the variation in the cohomology
  of the symplectic form of the reduced phase space}, Invent. Math. \textbf{69}
  (1982), no.~2, 259--268.

\bibitem[DH94]{dolg_hu}
Igor Dolgatchev and Yi~Hu, \emph{Variation of geometric invariant theory
  quotients}, eprint alg-geom/9402008, 1994.

\bibitem[DHVW85]{dhvw:i}
L.~Dixon, J.~A. Harvey, C.~Vafa, and E.~Witten, \emph{Strings on orbifolds},
  Nuclear Phys. B \textbf{261} (1985), no.~4, 678--686.

\bibitem[DK90]{don_kron:4mfds}
S.~K. Donaldson and P.~B. Kronheimer, \emph{The geometry of four-manifolds},
  Oxford University Press, 1990.

\bibitem[GH78]{gri_har:ag}
Phillip Griffiths and Joseph Harris, \emph{Principles of algebraic geometry},
  Wiley-Interscience (John Wiley \& Sons), New York, 1978.

\bibitem[GM87]{gold_mill:flat}
William~M. Goldman and John~J. Millson, \emph{Deformations of flat bundles over
  {K}{\"a}hler manifolds}, Geometry and Topology, Manifolds, Varieties and
  Knots (C.~McCrory and T.~Shifrin, eds.), Lecture Notes in Pure and Applied
  Mathematics, vol. 105, Marcel Dekker, New York-Basel, 1987, pp.~129--145.

\bibitem[GM88]{gold_mill:fundamental}
William~M. Goldman and John~J. Millson, \emph{The deformation theory of
  representations of fundamental groups of compact {K}{\"a}hler manifolds},
  Inst. Hautes Etudes Sci. Publ. Math. \textbf{67} (1988), 43--96.

\bibitem[GM90]{gold_mill:invariance}
William~M. Goldman and John~J. Millson, \emph{The homotopy invariance of the
  {K}uranishi space}, Illinois J. Math. \textbf{34} (1990), no.~2, 337--367.

\bibitem[Har77]{hart:ag}
Robin Hartshorne, \emph{Algebraic geometry}, Springer-Verlag, New
  York-Heidelberg, 1977.

\bibitem[It{\=o}94]{ito:trihedral}
Yukari It{\=o}, \emph{Crepant resolution of trihedral singularities},
  Proceedings of the Japan Academy \textbf{70} (1994), no.~Ser. A, no. 5,
  131--136.

\bibitem[Kro86]{kron:thesis}
Peter~B. Kronheimer, \emph{{ALE} gravitational instantons}, D. {P}hil. thesis,
  University of Oxford, 1986.

\bibitem[Kro89]{kron:ale}
P.~B. Kronheimer, \emph{The construction of {ALE} spaces as hyper-{K\"a}hler
  quotients}, J. Differential Geom. \textbf{29} (1989), no.~3, 665--683.

\bibitem[Kur62]{kura:complex}
M.~Kuranishi, \emph{On the locally complete families of complex analytic
  structures}, Ann. Math. (2) \textbf{75} (1962), 536--577.

\bibitem[Kur65]{kura:newproof}
M.~Kuranishi, \emph{New proof for the existence of locally complete families of
  complex structures}, Proc. {C}onf. {C}omplex {A}nalysis ({M}inneapolis,
  1964), Springer, Berlin, 1965, pp.~142--154.

\bibitem[Mar93]{mark:res_168}
D.~Markusevich, \emph{Resolution of ${H}_{168}$}, Preprint, Israel Institute of
  Technology, Haifa, 1993.

\bibitem[MOP87]{mar_ols_per}
D.~G. Markushevich, M.~A. Olshanetsky, and A.~M. Perelomov, \emph{Description
  of a class of superstring compactifications related to semisimple {L}ie
  algebras}, Comm. Math. Phys. \textbf{111} (1987), no.~2, 247--274.

\bibitem[Nad88]{nadel:quadratic}
A.M. Nadel, \emph{Singularities and {K}odaira dimension of the moduli space of
  {H}ermitian-{Y}ang-{M}ills connections}, Composito Math. \textbf{67} (1988),
  121--128.

\bibitem[Roa90]{roan:calabi-yau}
Shi~Shyr Roan, \emph{On {C}alabi-{Y}au orbifolds in weighted projective
  spaces}, Internat. J.~Math. \textbf{1} (1990), no.~2, 211--232.

\bibitem[Roa91]{roan:mirror_cy}
Shi~Shyr Roan, \emph{The mirror of {C}alabi-{Y}au orbifold}, Internat. J.~Math.
  \textbf{2} (1991), no.~4, 439--455.

\bibitem[Roa93]{roan:res_a5}
Shi~Shyr Roan, \emph{On $c_1=0$ resolution of quotient singularity}, Academia
  Sinica preprint R930826-2, 1993.

\bibitem[SI94]{sacha:thesis}
Alexander~V. Sardo~Infirri, \emph{Resolutions of orbifold singularities and
  representation moduli of {McK}ay quivers}, D. {P}hil. thesis, University of
  Oxford, 1994, Available as Pre-Print RIMS-984.

\bibitem[SI96a]{sacha:sl4}
Alexander~V. Sardo~Infirri, \emph{Crepant terminalisations and orbifold {E}uler
  number for {SL}$(4)$ singularities}, to appear, 1996.

\bibitem[SI96b]{sacha:flows}
Alexander~V. Sardo~Infirri, \emph{Resolutions of orbifold singularities and the
  transportation problem on the {M}c{K}ay quiver}, to appear, 1996.

\bibitem[Tha94]{thaddeus:git_flips}
M.~Thaddeus, \emph{Geometric invariant theory and flips}, eprint
  alg-geom/9405004, 1994.

\bibitem[Tia87]{tian:smoothness}
Gang Tian, \emph{Smoothness of the universal deformation space of compact
  {C}alabi-{Y}au manifolds and its {P}etersson-{W}eil metric}, Mathematical
  aspects of string theory (San Diego, Calif., 1986), Adv. Ser. Math. Phys.,
  vol.~1, World Sci. Publishing, Singapore, 1987, pp.~629--646.

\bibitem[Tod89]{todo:weil-peterson}
Andrey~N. Todorov, \emph{The {W}eil-{P}etersson geometry of the moduli space of
  {SU}$(n)$ $(n\geq 3)$ ({C}alabi-{Y}au) manifolds. {I}}, Comm. Math. Phys.
  \textbf{126} (1989), no.~2, 325--346.

\bibitem[UY86]{uhl_yau:hym}
K.~Uhlenbeck and S.~T. Yau, \emph{On the existence of
  {H}ermitian-{Y}ang-{M}ills connections in stable vector bundles}, Comm. Pure
  Appl. Math. \textbf{39} (1986), 257--.

\end{thebibliography}
\end{document}